\documentclass[11pt,twoside]{article}

\usepackage{epsfig,amsfonts,color}
\usepackage{amsmath}

\usepackage{amssymb, palatino, geometry,url}
\usepackage{algorithmic}
\usepackage[linesnumbered,ruled,noend,slide]{algorithm2e}

\usepackage[colorlinks=true,linkcolor=blue,citecolor=blue,urlcolor=blue]{hyperref}

\usepackage{algorithmic}
\usepackage{subfigure}
\usepackage[counterclockwise, figuresleft]{rotating}

\usepackage{fancyhdr}
\newcommand{\be}{\begin{equation}}
\newcommand{\ee}{\end{equation}}
\newcommand{\bea}{\begin{eqnarray}}
\newcommand{\eea}{\end{eqnarray}}

\newcommand{\bvec}{\left(\begin{array}{c}}
\newcommand{\evec}{\end{array}\right)}
\newcommand{\bsub}{\begin{subequations}}
\newcommand{\esub}{\end{subequations}}

\newcommand{\Lagr}{\mathcal{L}}

\DeclareMathOperator*{\argmin}{arg\,min}

\usepackage{lineno}
\usepackage[square,sort&compress,numbers]{natbib}
\makeatletter
\def\munderbar#1{\underline{\sbox\tw@{$#1$}\dp\tw@\z@\box\tw@}}
\makeatother

\begin{document}

\title{Scalable Preconditioning of Block-Structured\\  Linear Algebra Systems using ADMM}

\author{Jose S. Rodriguez${}^{\ddag}$, Carl D. Laird${}^{\dag}{}^{\ddag}$, and Victor M. Zavala${}^{\P}$\thanks{Corresponding Author: victor.zavala@wisc.edu}\\
 {\small ${}^\ddag$Davidson School Of Chemical Engineering }\\
 {\small Purdue University, 610 Purdue Mall, West Lafayette, IN 47907 }\\
 {\small ${}^\dag$Center for Computing Research}\\
  {\small Sandia National Laboratories,  Albuquerque, NM 87185}\\
  {\small ${}^{\P}$Department of Chemical and Biological Engineering}\\
 {\small \;University of Wisconsin-Madison, 1415 Engineering Dr, Madison, WI 53706, USA}}
 \date{}
\maketitle

\begin{abstract}
We study the solution of block-structured linear algebra systems arising in optimization by using iterative solution techniques. These systems are the core computational bottleneck of many problems of interest such as parameter estimation, optimal control, network optimization, and stochastic programming.  Our approach uses a Krylov solver (GMRES) that is preconditioned with an alternating method of multipliers (ADMM). We show that this ADMM-GMRES approach overcomes well-known scalability issues of Schur complement decomposition in problems that exhibit a high degree of coupling. The effectiveness of the  approach is demonstrated using linear systems that arise in stochastic optimal power flow problems and that contain up to 2 million total variables and 4,000 coupling variables. We find that ADMM-GMRES is nearly an order of magnitude faster than Schur complement decomposition.  Moreover, we demonstrate that the approach is robust to the selection of the augmented Lagrangian penalty parameter, which is a key advantage over the direct use of ADMM. 
\end{abstract}


\section{Introduction}\label{sec:intro}

The scalability of optimization solvers relies quite heavily on the solution of the underlying linear algebra systems. Advances in direct sparse linear algebra solvers have been instrumental in the widespread use of quadratic programming and nonlinear programming solvers such as {\tt Ipopt}, {\tt OOQP}, and {\tt Knitro} \citep{HSL, MUMPS1, MUMPS2}. Specialized direct solution techniques have also been developed to tackle large-scale and {\em block-structured} systems (using variants of Schur complement decomposition techniques) \citep{zavalalaird, WordDaniel2014Epso, pipsnlp, GondzioJacek2009Esip, GondzioJacek2003Pisf}.  Block structures appear in many important applications such as parameter estimation, stochastic programming, network optimization, and optimal control.  Schur decomposition techniques can also leverage {\em parallel computing architectures} and have enabled the solution of problems with millions to billions of variables and constraints.  Unfortunately, many applications of interest still remain inaccessible due to fundamental scalability limitations of Schur complement techniques. Specifically, Schur complement decomposition does not scale well in problems that exhibit high degrees of {\em block coupling}.  This is because high degrees of coupling require assembling and factorizing large Schur complement matrices (which are often highly dense). 

Iterative solution techniques \cite{QuarteroniAlfio2007Nm, MaC.F.2015TcUm, ForsgrenAnders2007ISoA, ElmanHowardC.1994IaPU, BenziMichele2006Oteo, BenziMichele2005Nsos} and associated preconditioning strategies  \cite{CaoYankai2016Cpfs, MoralesJoseLuis2000APbL, GolubGeneH.2003OSBI, RustenTorgeir1992APIM, GillPhilipE.1992PfIS, WalterZulehner2002Aoim} have been proposed to address fundamental scalability issues of direct linear algebra strategies. In the context of block-structured problems, attempts have been made to solve the Schur complement system by using iterative solution techniques (to avoid assembling and factorizing the Schur complement). Preconditioners for Schur complements arising in special problem classes such as multi-commodity networks and stochastic programs have been developed \cite{CaoYankai2016Cpfs}. Unfortunately, preconditioning strategies for general problem classes are still lacking.  Another important issue that arises in this context is that the implementation of advanced linear algebra strategies is non-trivial (e.g., it requires intrusive modifications of optimization solvers). 

Along a separate line of research, significant advances have been made in the development of {\em problem-level} decomposition techniques such as the alternating direction method of multipliers (ADMM) and Lagrangian dual decomposition \cite{HanDeren2013LLCo, GuoKe2017CoAf, HongMingyi2017Otlc, HanDeren2013LLCo, GoldfarbDonald2012FMAf, HeBingsheng2012OtOC, RodriguezJoseS.2018BAin}. Such approaches are flexible and rather easy to implement but suffer from slow convergence.  Recently, it has been proposed to use ADMM as a preconditioner for Krylov-based iterative solvers such as GMRES \cite{ZhangRichard2016PIiA, ZhangRichardY.2018GAfQ}. In this work, we provide a detailed derivation of this ADMM-GMRES approach and test its performance in the context of block-structured linear algebra systems.  We demonstrate that this approach overcomes the scalability issues of Schur complement decomposition. We also demonstrate that this approach is significantly more effective than using ADMM directly. Our tests are facilitated by the use of {\tt PyNumero}, a recently-developed Python framework that enables the implementation and benchmarking of optimization algorithms. We use the proposed framework to tackle problems with hundreds of thousands to millions of variables generated from standard benchmark sets and power grid applications.  

The paper is structured as follows. In Section \ref{sec:prelims} we define the problem of interest and provide preliminary information on the use of Schur complement decomposition and ADMM approaches. In Section \ref{sec:derivation} we provide a detailed derivation of the ADMM-GMRES approach and in Section \ref{sec:benchmarks} we provide benchmark results. 

\section{Preliminaries}\label{sec:prelims}

We study the solution of block-structured quadratic programs (QP) of the form:
\begin{subequations}
\begin{align}
\underset{x_i, z}{\text{min}}
 \;\; & \;\;   \sum_{i \in \mathcal{P}} \frac{1}{2}x_i^T D_i x_i + c_i^T x_i& \label{eq:QB2-a} \\
\text{s.t.} \;\; &\;\; J_i x_i  = b_i,\; \;\;\;\;\;\; \quad (\lambda_i)\quad i \in \mathcal{P} \label{eq:QB2-b}\\
&\;\; A_{i}x_i + B_{i}z = 0,\; (y_i)\quad i \in \mathcal{P} \label{eq:QB2-c}.
\end{align}
\label{eq:block-structured-qp}
\end{subequations}
\noindent Here, $\mathcal{P} := \{1,\dots,P\}$ is a set of block variable partitions. Each partition contains a vector of primal variables $ x_i \in \mathbb{R}^{n_{x_i}}$ and the vector $ z \in \mathbb{R}^{n_z}$ contains the primal variables that the couple partitions. The total number of primal variables is $n := n_z + \sum_{i\in\mathcal{P}}n_{x_i}$. Equation \eqref{eq:QB2-b} are the partition constraints with their respective dual variables $\lambda_i\in \mathbb{R}^{m_i}$. Equation \eqref{eq:QB2-c} are the constraints that {\em link} partitions across set $\mathcal{P}$ and have associated dual variables $y_i\in\mathbb{R}^{l_i}$. We assume that the partition matrices $J_i\in\mathbb{R}^{m_i\times n_i}$ have full rank and that the right-hand-side coefficients $b_i\in \mathbb{R}^{m_i}$ are in the column space of $J_i$. The total number of partition constraints is $m :=\sum_{i\in\mathcal{P}}m_i$. We refer to $A_i\in \mathbb{R}^{n_z\times n_i}$ and $B_i\in \mathbb{R}^{n_z\times n_z}$ as linking matrices and we assume them to have full rank. The total number of linking constraints is $l:=\sum_{i\in\mathcal{P}}l_i$.  The QP under study is the main computational kernel behind nonlinear programming strategies because it is used for computing primal-dual search steps. 

We make the blanket assumption that the block-structured QP is strongly convex and that the combined Jacobian matrix (obtained by assembling partition and coupling constraints) has full rank. Strong convexity can be obtained by ensuring that all block Hessian matrices $D_i$ are positive definite. Strong convexity and full-rank conditions guarantee that the primal-dual solution of the QP exists and is unique. Moreover, these assumptions guarantee that the QP solution is a unique minimizer and that this can be found by solving the first-order stationarity conditions.  Additional assumptions will also be needed on the nature of the building blocks of the QP (associated with each partition). Such assumptions are needed to ensure that proposed decomposition schemes are well-defined and will be stated as we proceed (in order to maintain clarity in the presentation). 

\noindent The Lagrange function of \eqref{eq:block-structured-qp} can be expressed as:
\begin{align}
\label{eqn:auglag}
\Lagr(x,z,\lambda, y) =  \sum_{i \in \mathcal{P}} \frac{1}{2}x_i^T D_i x_i + c_i^T x_i  + y_i^T(A_ix_i+B_i z)+ \lambda_i^{T}(J_ix_i-b_i),
\end{align}
\noindent where $x:=(x_1, \cdots, x_P)$, $\lambda:=(\lambda_1, \cdots, \lambda_P)$ and $y:=(y_1, \cdots, y_P)$. The first-order optimality conditions of \eqref{eq:block-structured-qp} are given by:
\begin{eqnarray}
\nabla_{x_i}\mathcal{L}&=0 =&D_ix_i + c_i + J_i^T\lambda_i + A_i^{T}y_i,\; i \in \mathcal{P}\nonumber\\
\nabla_{\lambda_i}\mathcal{L}&=0=&J_ix_i-b_i,\; i \in \mathcal{P} \nonumber\\
\nabla_{y_i}\mathcal{L}&=0=& A_ix_i-B_iz,\; i \in \mathcal{P}\nonumber\\
\nabla_z\mathcal{L}&=0=&\sum_{i \in \mathcal{P}} B_i^Ty_i.\label{eq:ocp-qp}
\end{eqnarray}
\noindent These conditions form a {\em block-structured} linear system of the form shown in \eqref{eq:kkt-blocks}. For the sake of compactness and ease of notation, we rewrite \eqref{eq:kkt-blocks} as:
\begin{equation}
\label{eq:kkt-compact}
\underbrace{\begin{bmatrix}
K &   &  A^{T}   \\
 &                & B^{T}                \\
A  &   B      &                 \\   
\end{bmatrix}}_{H}
\underbrace{\begin{bmatrix}
\upsilon\\
z\\
y\\
\end{bmatrix}}_u
=
\underbrace{\begin{bmatrix}
\gamma \\
0 \\
0 \\
\end{bmatrix}}_r.
\end{equation}
\noindent where $\upsilon=(\upsilon_1, \cdots, \upsilon_P)$, $\upsilon_i = (x_i, \lambda_i)$,  $\gamma=(\gamma_1, \cdots, \gamma_P)$, $\gamma_i = (-c_i, b_i)$, $u=(v,z,y)$, $r=(\gamma,0,0)$, and $K = \text{blkdiag}\{K_1, K_2, \cdots, K_P\}$. We also have $A = \text{blkdiag}\{\tilde{A}_1, \tilde{A}_2, \cdots, \tilde{A}_P\}$, $B = \text{rowstack}\{B_1, B_2, \cdots, B_P\}$ with:

\begin{equation}
K_i=\begin{bmatrix}
D_i &  J_i^{T} \\
J_i  &  \\
\end{bmatrix}            
\qquad \tilde{A}_i=\begin{bmatrix}
A_i &  \quad 0
\end{bmatrix}          ,\; i \in \mathcal{P}.
\end{equation}  

\begin{equation}
\label{eq:kkt-blocks}
\begin{bmatrix}
D_1 &  J_1^{T} &            &           &               &          & A_1^{T}  &             & \\
J_1  &               &            &           &               &         &                &              & \\
       &               &   \ddots &          &               &         &                &   \ddots &  \\
       &               &             &    D_P &  J_P^{T} &         &                &              & A_P^{T} \\
      &               &             &      J_P &               &          &               &              &  \\
      &               &             &           &                &          &  B_1^{T} &  \cdots   & B_P^{T}\\
A_1 &              &            &            &                &  B_1  &               &               & \\
       &              &  \ddots &           &                & \vdots &             &                &\\
       &              &            &  A_P    &                &  B_P    &             &                &\\
\end{bmatrix}
\begin{bmatrix}
x_1\\
\lambda_1\\
\vdots\\
x_P\\
\lambda_P\\
z\\
y_1\\
\vdots\\
y_P
\end{bmatrix}=
\begin{bmatrix}
-c_1\\
b_1\\
\vdots\\
-c_P\\
b_P\\
0\\
0\\
\vdots\\
0
\end{bmatrix}
\end{equation}

\subsection{Solution using Schur Decomposition}

One can solve large instances of the block-structured QP by using a Schur-complement decomposition method (we refer to this approach simply as Schur decomposition) \cite{schurbook}. This approach decomposes \eqref{eq:block-structured-qp} by using block Gaussian elimination on a permuted version of the linear system \eqref{eq:kkt-blocks}. The permuted system has the structure:

\begin{equation}
\label{eq:kkt-blocks-for-sc}
\begin{bmatrix}
D_1 &  J_1^{T} &  A_1^{T} &               &            &          &                &             & \\
J_1  &               &                &               &               &         &                &              & \\
A_1 &               &                &               &               &         &                &       &  B_1\\
       &              &                  &   \ddots &               &           &               &               &  \vdots\\
       &               &                 &              & D_P       &  J_P^{T} &  A_P^{T}       &                &  \\
       &               &                 &              &  J_P      &               &              &         \\
       &               &                 &              &  A_P      &                &          &               &  B_P\\
       &              &    B_1^{T}  & \cdots   &                &            &    B_P^{T}     &                &\\
\end{bmatrix}
\begin{bmatrix}
x_1\\
\lambda_1\\
y_1\\
\vdots\\
x_P\\
\lambda_P\\
y_P\\
z
\end{bmatrix}=
\begin{bmatrix}
-c_1\\
b_1\\
0\\
\vdots\\
-c_P\\
b_P\\
0\\
0
\end{bmatrix}
\end{equation}

\noindent This system can be expressed in compact form as:

\begin{equation}
\label{eq:kkt-compact-sc}
\begin{bmatrix}
K_s & B_s    \\
B_s^{T}  &                   \\   
\end{bmatrix}
\begin{bmatrix}
\upsilon_s\\
z\\
\end{bmatrix}
=
\begin{bmatrix}
\gamma_s \\
0 \\
\end{bmatrix}
\end{equation}

\noindent where  $\upsilon_s=(\upsilon_{s_1}, \cdots, \upsilon_{s_P})$, $\upsilon_{s_i} = (x_i, \lambda_i, y_i)$,  $\gamma_s=(\gamma_{s_1}, \cdots, \gamma_{s_P})$, and $\gamma_i = (-c_i, b_i, 0)$. We also have $K_s = \text{blkdiag}\{K_{s_1}, K_{s_2}, \cdots, K_{s_P}\}$ and $B_s = \text{rowstack}\{B_{s_1}, B_{s_2}, \cdots, B_{s_P}\}$ with:
\begin{equation}
K_{s_i}=\begin{bmatrix}
D_i &  J_i^{T} & A^{T} \\
J_i  &  &\\
A &
\end{bmatrix} \qquad B_{s_i}=\begin{bmatrix}
0&  0& B_i^{T} \\
\end{bmatrix}  ^{T} ,\; i \in \mathcal{P}.
\end{equation}

\noindent Because we have assumed that $D_i$ is positive definite and the combined constraint Jacobian $(J_i,A)$ has full rank, we have that $K_{s_i}$ is nonsingular. This also implies $K_s$ is nonsingular and,  as a result, we can form the Schur complement system:
\begin{equation}
\label{eq:schur-system}
(B_s^{T}K_s^{-1}B_s) z = B_s^{T}K_s^{-1}\gamma_s.
\end{equation}
\noindent We refer to the coefficient matrix of \eqref{eq:schur-system} as the Schur complement. Since $K_s$ is block-diagonal, one can assemble the Schur complement by factorizing the blocks $K_{s_i}$ independently. By using this assembled Schur complement, one can then factorize the Schur complement matrix and solve system \eqref{eq:schur-system} to find a solution for the coupling variables $z$. Having $z$, one then proceeds to find solutions for the partition variables $\upsilon_s$ by solving the following system:
\begin{equation}
\label{eq:schur-backsolves}
K_{s} \upsilon_{s} = \gamma_{s} - B_s z.
\end{equation}
\noindent Here, again, one can solve for each element $\upsilon_{s_i}$ independently because $K_s$ is block-diagonal. The  Schur decomposition method is summarized in Algorithm \ref{alg:sc-pseudo}. We refer the reader to \citep{Kang2014} for details on the implementation of Schur decomposition approaches. 

\begin{center}
\begin{algorithm}[ht!]
\label{alg:sc-pseudo}
\caption{Schur Decomposition for Block-Structured QP}
Let $S=0$ and $r_{sc}=0$\\
  \small
	Factorize $K_s$ matrix :   \\
		\ForEach{ $i \in \mathcal{P}$}{
			Factorize $K_{s_i}$\\
}
  Form Schur complement system: \\		
 \ForEach{ $i \in \mathcal{P}$}{
 $S = S + B_i^{T}K_{s_i}^{-1}B_i$\\
 $r_{sc} = r_{sc} + B_i^{T}K_{s_i}^{-1}\gamma_{s_i}$\\
  }
 Factorize $S$ and compute coupling variables by solving: \\
 $\quad Sz=r_{sc}$\\
 Compute partition variables: \\
  \ForEach{ $i \in \mathcal{P}$}{
  $\quad K_{s_i} \upsilon_{s_i} = \gamma_{s_i} - B_iz$\\
  }
\end{algorithm}
\end{center}

Schur decomposition is a flexible approach that enables the solution of problems with many block partitions. A fundamental limitation of this approach, however, is that one needs to assemble and factorize the Schur complement matrix $S$ (which is often highly dense). As a result, Schur decomposition does not scale well with the number of coupling variables $z$. Iterative approaches can in principle be used to solve the Schur system but effective preconditioning strategies do not currenlty exist for general block-structured systems. 

\subsection{Solution using ADMM}

The block-structured QP can also be decomposed and solved by using ADMM. This approach seeks to minimize the augmented Lagrangian function:
\begin{align}
\Lagr_\rho(x,z,\lambda, y)=\sum_{i\in \mathcal{P}}x_i^TD_ix_i + c_i^Tx_i  +\left(A_ix_i  + B_iz\right)^{T}y_i + \lambda_i^{T}(J_ix_i-b_i)+  \frac{\rho}{2}\|A_ix_i  + B_iz\|^{2}
\end{align}
by performing alternating minimization with respect to the block variables $(x_i,y_i)$ and the coupling variables $z$. A standard implementation of ADMM for solving the QP of interest is presented in Algorithm \ref{alg:admm-qp}. 
\begin{center}
\begin{algorithm}[ht!]
\label{alg:admm-qp}
\caption{ADMM for Block-Structured QP}
  \small
  \KwIn{Starting point $u^{0} = (\upsilon^{0}, y^{0}$, $z^{0})$, maximum number of iterations $N_{\textrm{ADMM}}$, penalty parameter $\rho > 0$, and convergence tolerance $\epsilon > 0$}
   \For{ $k = 0, 1, 2, \ldots, N$}{
	Update partition variables:   \\
		\ForEach{ $i \in \mathcal{P}$}{
$\quad x_{i}^{k+1} = \argmin\limits_{ x_i \in \mathcal{X}_i}  x_i^TD_ix_i + c_i^Tx_i  +\left(A_ix_i  + B_iz^{k}\right)^{T}y_i^{k} + \frac{\rho}{2}\|A_ix_i  + B_iz^{k}\|^{2}$\\
}
  Update coupling variables: \\		
$\quad z^{k+1}= \argmin\limits_{z}  \left(A_ix_i^{k+1}  + B_iz\right)^{T}y_i^{k} + \frac{\rho}{2}\|A_ix_i^{k+1}  + B_iz\|^{2}$\\
 Update dual variables:\\
 \ForEach{ $i \in \mathcal{P}$}{
 $\quad y_i^{k+1} = y_i^{k} + \rho\left(A_ix_i^{k+1}  + B_iz^{k+1}\right)$\\
  }
 \If{$\|y^{k+1}-y^{k}\| \leq \epsilon$ \text{and} $\|\rho A^{T}B\cdot(z^{k+1}-z^{k})\| \leq \epsilon$}{
 	\textbf{stop}\\
 }
}  
\KwOut{$u$} 
\end{algorithm}
\end{center} 
In the above algorithm, $\mathcal{X}_i=\{x\,|\,J_ix-b_i=0\}$ is used to denote the feasible set of each partition (the inner block constraints are satisfied exactly). The ADMM algorithm can be implemented by solving the first-order conditions of each subproblem directly. This is because each block subproblem is strongly convex and the block Jacobian has full rank. This approach is sketched in Algorithm \ref{alg:admm}. 

\begin{center}
\begin{algorithm}[H]
\label{alg:admm}
\caption{{\tt ADMM}($u^{0}$, $N$, $\rho$)}
  \small
  \KwIn{ starting point $u^{0} = (\upsilon^{0}, y^{0}$, $z^{0})$, maximum number of iterations $N_{\textrm{ADMM}}$, penalty parameter $\rho > 0$,  and convergence tolerance $\epsilon > 0$ }
	Factorize $K_{\rho}$ and $B^{T}B$ matrix:\\
	\ForEach{ $i \in \mathcal{P}$}{
	Factorize partition matrices $K_{\rho_i}$ and $B_i^{T}B_i$\\
}  
   \For{ $k = 0, 1, 2, \ldots, N$}{
	Update partition variables :   \\
		\ForEach{ $i \in \mathcal{P}$}{
$\quad K_{\rho_i}\upsilon^{k+1} = -\left(\gamma_i + \begin{bmatrix}\rho A_i^{T}B_i z^{k} \\ 0\end{bmatrix} + \begin{bmatrix}A_i^{T}y_i^{k} \\ 0\end{bmatrix}\right)$\\
}
  Update coupling variables: \\		
$\quad z^{k+1}= -[B^TB]^{-1}\left(B^{T}A\upsilon^{k+1}+\frac{1}{\rho}B^{T}y^{k}\right)$\\
 Update dual variables:\\
 \ForEach{ $i \in \mathcal{P}$}{
 $\quad y_i^{k+1} = y_i^{k} + \rho\left(\tilde{A}_i\upsilon_i^{k+1}  + B_iz^{k+1}\right)$\\
  }
 \If{$\|y^{k+1}-y^{k}\| \leq \epsilon$ \text{and} $\|\rho A^{T}B\cdot(z^{k+1}-z^{k})\| \leq \epsilon$}{
 	\textbf{stop}
 }
}  
\KwOut{$u$} 
\end{algorithm}
\end{center}

In the above algorithm we have that $K_{\rho} = \text{blkdiag}\{K_{\rho_1}, K_{\rho_2}, \cdots, K_{\rho_P}\}$ with
\begin{equation}
K_{\rho_i}=\begin{bmatrix}
D_i + \rho A_i^TA_i&  J_i^{T} \\
J_i  & 0 \\
\end{bmatrix}.            
\end{equation}

We note that the update of the coupling variables still requires forming and factorizing the matrix $B^TB$ (but this can be done in blocks by forming and factorizing $B_i^TB_i$ individually). Moreover, this operation only needs to be performed once. Note also that $B^TB$ is invertible since $B$ has full rank. As a result, the update step for the coupling variables in ADMM is cheaper than that of Schur decomposition and can thus overcome the main computational bottleneck of the latter. Unfortunately, it is well-known that ADMM exhibits slow convergence and thus the ability to perform fast operations might be shadowed by the need to perform many iterations. 

\section{Solution using ADMM-GMRES}\label{sec:derivation}

The key observation that motivates our work is that ADMM can be used as a {\em preconditioner} for iterative linear algebra techniques such as GMRES \cite{ZhangRichard2016PIiA, ZhangRichardY.2018GAfQ}. To derive the ADMM preconditioning strategy, we consider the {\em regularized} QP \eqref{eq:block-structured-qp}:

\begin{subequations}
\begin{align}
\underset{x_i, z}{\text{min}}
 \;\; & \;\;   \sum_{i \in \mathcal{P}} \frac{1}{2}x_i^T D_i x_i + c_i^T x_i + \frac{\rho}{2}\|Ax_i+B_iz\|^{2}& \label{eq:AUGQB2-a} \\
\text{s.t.} \;\; &\;\; J_i x_i  = b_i,\; \;\;\;\;\;\; \quad (\lambda_i)\quad i \in \mathcal{P} \label{eq:AUGQB2-b}\\
&\;\; A_{i}x_i + B_{i}z = 0,\; (y_i)\quad i \in \mathcal{P} \label{eq:AUGQB2-c}.
\end{align}
\label{eq:block-structured-augmented-qp}
\end{subequations}
\noindent The solution of this problem is also a solution of \eqref{eq:block-structured-qp} (since the penalization term vanishes at the solution). The optimality conditions of the regularized QP are given by:

\begin{equation}
\label{eq:kkt-aug-compact}
\underbrace{\begin{bmatrix}
K_{\rho} &  \rho A^{T}B &  A^{T}   \\
 \rho B^TA &       \rho B^{T}B         & B^{T}                \\
A  &   B      &                 \\   
\end{bmatrix}}_{H_\rho}
\underbrace{\begin{bmatrix}
\upsilon\\
z\\
y\\
\end{bmatrix}}_{u}
=
\underbrace{\begin{bmatrix}
\gamma \\
0 \\
0 \\
\end{bmatrix}}_r.
\end{equation}

We refer to \eqref{eq:kkt-aug-compact} as the \textit{KKT system} and to $H_{\rho}$ as the \textit{KKT matrix}.   ADMM can be interpreted as a Gauss-Seidel (alternating) minimization of the block and coupling variables and the dual variables \citep{boyd2011distributed, eckstein1992douglas, RodriguezJoseS.2018BAin}. This induces a splitting operator $H_{\rho}=M_{\rho}-N_{\rho}$ satisfying:

\begin{equation}
\label{eq:splitting}
\underbrace{\begin{bmatrix}
K_{\rho} &  \rho A^{T}B &  A^{T}   \\
 \rho B^TA&       \rho B^{T}B         & B^{T}                \\
A  &   B      &                 \\   
\end{bmatrix}}_{H_\rho}=
\underbrace{\begin{bmatrix}
K_{\rho} &   & \\
\rho B^TA &       \rho B^TB   &               \\
A  &   B      &       -\frac{1}{\rho}I           \\   
\end{bmatrix}}_{M_{\rho}}-
\underbrace{\begin{bmatrix}
\; &  -\rho A^TB &  -A^{T}   \\
\; &                   & -B^{T}                \\
\; &                   &  -\frac{1}{\rho}I          \\   
\end{bmatrix}}_{N_{\rho}}
\end{equation}
\noindent Applying splitting \eqref{eq:splitting} to \eqref{eq:kkt-aug-compact} gives the operator:
\begin{equation}
\label{eq:fixed-point}
T_\rho(u) := G_{\rho}u + f_{\rho}
\end{equation}
\noindent where $G_{\rho}{=}M_{\rho}^{-1}N_{\rho}$ and $f_{\rho} {=} M_{\rho}^{-1}r$. Note that any $u$ satisfying the fixed-point $T_\rho(u){=}u$ also satisfies $(I-G_{\rho})u {=} f_{\rho}$ and is a solution of the preconditioned KKT system:
\begin{equation}
\label{eq:preconditioned-sys}
M_{\rho}^{-1}H_{\rho}u = M_{\rho}^{-1}r.
\end{equation}
\noindent This follows from $M_{\rho}^{-1}H_{\rho} = M_{\rho}^{-1}\left(M_{\rho}-N_{\rho}\right) = I-G_{\rho}$. This motivates the development of a Richardson recursion of the form $u^{k+1} {=} G_{\rho}u^{k} + f_{\rho}$, which converges to a $u$ satisfying $(I-G_{\rho})u {=} f_{\rho}$ and $M_{\rho}^{-1}H_{\rho}u {=} M_{\rho}^{-1}r$ (provided that the eigenvalues of $G_{\rho}$ are inside the unit circle). In Appendix \ref{app:admm-richardson} we show that the operator $T_\rho(u)$ can be computed by performing one ADMM iteration (using $u$ as starting point).  In other words, we have that $T_\rho(u){=}${\tt ADMM}($u$, $N{=}1$, $\rho$). This also implies that the Richardson recursion can be written as $u^{k+1}{=}${\tt ADMM}($u^k$, $N{=}1$, $\rho$). Consequently, the Richardson recursion (and thus ADMM) are consistent preconditioner choices. 

The key idea behind ADMM-GMRES is to solve the system $M_{\rho}^{-1}H_{\rho}u {=} M_{\rho}^{-1}r$ by using the Krylov solver GMRES. This is equivalent to solving $(I-G_{\rho})u {=} f_{\rho}$. The right-hand side of this system can be computed as $ f_{\rho}=T_\rho(0)$. GMRES requires the computation of matrix-vector products with the preconditioned coefficient matrix of the form $M_{\rho}^{-1}H_{\rho}h{=}(I-G_\rho)h$. This can be done by using the operator \eqref{eq:fixed-point} as:
\begin{align}
M_{\rho}^{-1}H_{\rho}h &= h - \left[T_\rho(h) - T_\rho(0)\right],
\end{align}
This follows from the observation that:
\begin{align}\label{eq:matvec}
M_{\rho}^{-1}H_{\rho}h&=h - \left[T_\rho(h) - T_\rho(0)\right]\nonumber\\
&= h - \left[G_{\rho}h+f_{\rho} -f_{\rho}\right]\nonumber\\
&= h - G_{\rho}h \nonumber\\
&=(I-G_{\rho})h.
\end{align}
From \eqref{eq:matvec} we note that asking the ADMM oracle {\tt ADMM}($h$, $N$, $\rho$) to iterate until reaching convergence will deliver $T_\rho(h){=}h$ satisfying $M_{\rho}^{-1}H_{\rho}h{=}f_\rho$. In such a case, the ADMM preconditioner is perfect (since it solves the actual preconditioned KKT system). Consequently, the quality of the ADMM preconditioner will improve as one increases $N$. For details on the properties of the preconditioner, we refer the reader to \citep{ZhangRichardY.2018GAfQ, ZhangRichard2016PIiA}. The ADMM-GMRES strategy is summarized in Algorithm \ref{alg:admm-gmres}.

\begin{center}
\begin{algorithm}[H]
\label{alg:admm-gmres}
\caption{{\tt ADMM\_GMRES}($N_{\text{GMRES}}$, ${N}_{\text{ADMM}}$, $\rho$)}
  \small
  \KwIn{maximum number of GMRES iterations $N_{\text{GMRES}}$, maximum number of ADMM iterations ${N}_{\text{ADMM}}$, penalty parameter $\rho > 0$, and tolerance $\epsilon > 0$}
  Compute right-hand-side vector:\\
  $\quad f_\rho = $ {\tt ADMM}(0, ${N}_{\text{ADMM}}$, $\rho$)\\
  Call GMRES solver\footnote{Matrix-vector products are computed as $(I-G_{\rho})h {=} h-(T_\rho(h) - f_\rho)$, where $T_\rho(h){=}${\tt ADMM}($h$, ${N}_{\text{ADMM}}{=}1$, $\rho$).}:\\
  $u$ = {\tt GMRES}$(I-G_{\rho}, f_\rho, N_{\text{GMRES}}, \epsilon)$\\
  \KwOut{$u$} 
\end{algorithm}
\end{center}

\section{Numerical Results}\label{sec:results}

In this section we discuss the implementation of ADMM-GMRES and present results for different benchmark problems. All numerical experiments were performed using {\tt PyNumero}, which is an open-source framework written in Python and C++ that combines modeling capabilities of the algebraic modeling language {\tt Pyomo} \citep{hart2017pyomo} with efficient libraries like the {\tt AMPL} solver library \citep{ampl}, the Harwell Subroutine Library (HSL), and NumPy/SciPy \citep{scipyref}. It uses object-oriented principles that facilitate the implementation of algorithms and problem formulations that exploit block-structures via polymorphism and inheritance. All these features facilitate the implementation of ADMM, Schur decomposition, and ADMM-GMRES. The optimization models were implemented in {\tt Pyomo}/{\tt PyNumero} and all linear algebra operations were performed in compiled code. Within {\tt PyNumero}, we used an {\tt MA27} interface to perform all direct linear algebra operations. We use the {\tt GMRES} implementation available in {\tt Scipy} to perform all iterative linear algebra operations. Iterative linear algebra routines available in {\tt KRYPY} \citep{GaulAndre2015PRKs} were also used to validate results. To implement the power grid models we used {\tt EGRET} \footnote{\url{https://github.com/grid-parity-exchange/Egret}}, a {\em Pyomo}-based package that facilitates the formulation if optimization problems that arise in power systems. The convergence criterion for GMRES and ADMM requires that the norm of the KKT system residual $H u-r$ is smaller than $\epsilon=1\cdot 10^{-8}$.  If the convergence criterion is not satisfied after 2,000 iterations, the algorithm was aborted and we report the final residual achieved. The linear solver {\tt MA27} was used with a pivoting tolerance of $1\cdot10^{-8}$

\subsection{Standard Benchmark Problems}\label{sec:benchmarks}

We first conducted tests with randomly generated instances to study qualitatively the performance of ADMM-GMRES on block-structured optimization problems. This section focuses on two-stage stochastic programs of the form:
\begin{subequations}
\begin{align}
\underset{x_i, z}{\text{min}}
 \;\; & \;\;   \sum_{i \in \mathcal{P}} \frac{1}{2}x_i^T D x_i + c^T x_i& \label{eq:stoch-QB2-a} \\
\text{s.t.} \;\; &\;\; J x_i  = b_i,\; \;\;\;\;\;\; \quad (\lambda_i)\quad i \in \mathcal{P} \label{eq:stoch-QB2-b}\\
&\;\; A_{i}x_i - z = 0,\;\;\;\;\; (y_i) \quad i \in \mathcal{P} \label{eq:stoch-QB2-c}
\end{align}
\label{eq:stochastic-qp}
\end{subequations}
\noindent where $\mathcal{P}$ is the scenario set, $x_i$ are the second-stage (recourse) variables, and $z$ are first-stage (coupling) variables. We defined a nominal vector $b$ and create scenarios with right-hand-side vector $b_i$ using the nominal vector $b$ as the mean and a standard deviation $\sigma{=}0.5b$. We first demonstrate the scalability of Algorithm \ref{alg:admm-gmres} when solving instances of problem \eqref{eq:stochastic-qp} with high dimensionality in the coupling variables $z$. The stochastic problem was constructed in the following manner: $D$ was set to a $4{,}800\times4{,}800$ block diagonal matrix with $16$ dense symmetric blocks. Each block was generated following Algorithm 14 from \citep{ZhangRichardY.2018GAfQ} with log-standard-deviation $s{=}0.5$ (see \citep{ZhangRichardY.2018GAfQ} for details). The random matrix $J$ has dimensions of $100\times4{,}800$. The number of scenarios was set to $50$, giving an initial problem with $240{,}000$ variables and $5{,}000$ constraints. To investigate the scalability of Algorithm \ref{alg:admm-gmres}, the number of complicating variables was varied from $100$ to $4{,}000$. Note that, as $n_z$ increases, the number of constraints of \eqref{eq:stochastic-qp} increases as $50n_z$. The largest problem solved contained $244{,}000$ total variables and $205{,}000$ total constraints.

\begin{figure}[!htb]
  \begin{center}
   \includegraphics[width=0.7\textwidth]{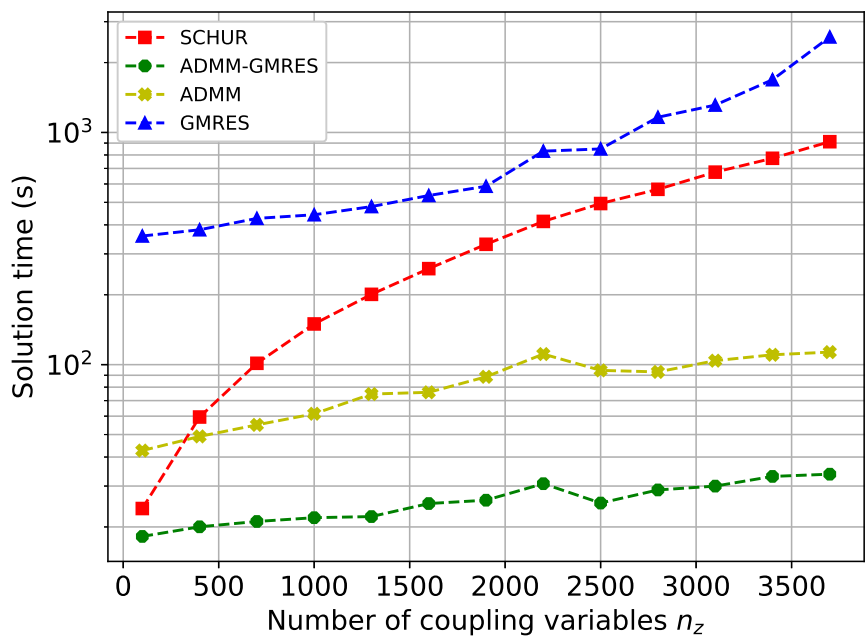}
 \caption{Scalability analysis of Schur decomposition, GMRES (without preconditioner), ADMM, and ADMM-GMRES.}\label{fig:scalability-results} 
 \end{center}
\end{figure}

We solved \eqref{eq:stochastic-qp} using four different strategies. Figure \ref{fig:scalability-results} summarizes the results obtained with Schur decomposition, GMRES (without preconditioner), ADMM, and ADMM-GMRES. These results confirm the observations of Section \ref{sec:prelims}. Specifically, Schur decomposition does not scale well as $n_z$ increases. The main reason for this behavior is that, as $n_z$ increases, the number of operations required to form the Schur-complement increase. In addition, because the Schur-complement matrix is a dense $n_z\times n_z$ matrix, the factorization time increases cubically as $n_z$ increases. We observe that GMRES takes the longest time to solve the problem. For ADMM and ADMM-GMRES we see more favorable scalability as $n_z$ increases. Specifically, we see that ADMM-GMRES converges faster and that the savings increase as the number of coupling variables increases.  We also note that ADMM-GMRES mimics the trend in performance of ADMM but is significantly faster.

\begin{figure}[!htb]
  \begin{center}
   \includegraphics[width=0.7\textwidth]{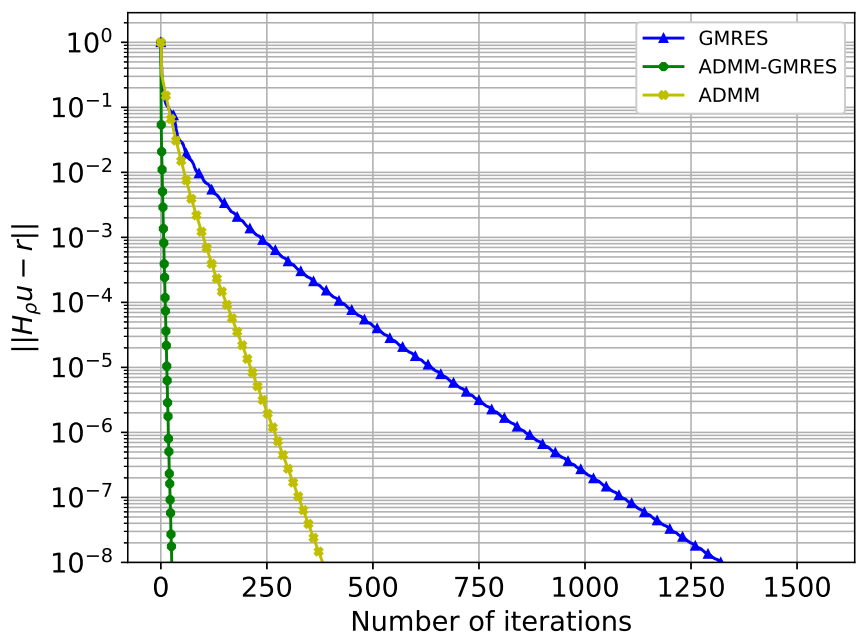}
 \caption{Evolution of residuals for GMRES (without preconditioner), ADMM, and ADMM-GMRES.}\label{fig:residuals}  
 \end{center}
\end{figure}

Figure \ref{fig:residuals} shows the residuals for the iterative approaches for a problem with $1{,}000$ complicating variables. We can see that all methods exhibit linear convergence but that ADMM-GMRES outperforms both ADMM and GMRES (unpreconditioned). Notably, ADMM-GMRES converges in just $35$ iterations while ADMM and GMRES require over $300$ iterations and $1{,}000$ iterations, respectively.  Moreover, we note that ADMM-GMRES can reach high accuracy levels (of $1\cdot 10^{-8}$), which is a desirable feature of iterative solution strategies.

An important drawback of ADMM is the need to tune the penalty parameter $\rho$. The work in \citep{GhadimiEuhanna2015OPSf} shows that an optimal value for $\rho$ can be chosen based on the smallest and largest eigenvalues of the matrix  $A^{T}D^{-1}A$. In principle, this selection of $\rho$ is optimal for ADMM-GMRES as well. However, for large-scale structured problems such as the ones considered here, computing the eigenvalues of $A^{T}D^{-1}A$ is expensive. Heuristic approaches have also been proposed to select $\rho$ at every ADMM iteration with the objective of accelerating convergence \citep{BrendtWohlberg2017APPS}. Unfortunately these heuristics do not provide guarantees and might incur additional overhead. In particular, for the QP problems considered here, varying $\rho$ at every ADMM iteration will require forming and factorizing the $K_{\rho_i}$ repetitively. Interestingly, we now proceed to demonstrate that ADMM-GMRES is fairly insensitive to the choice of $\rho$. 

\begin{figure}[!h]
  \centering
  \subfigure[Number of iterations]{{\includegraphics[width=3.25in]{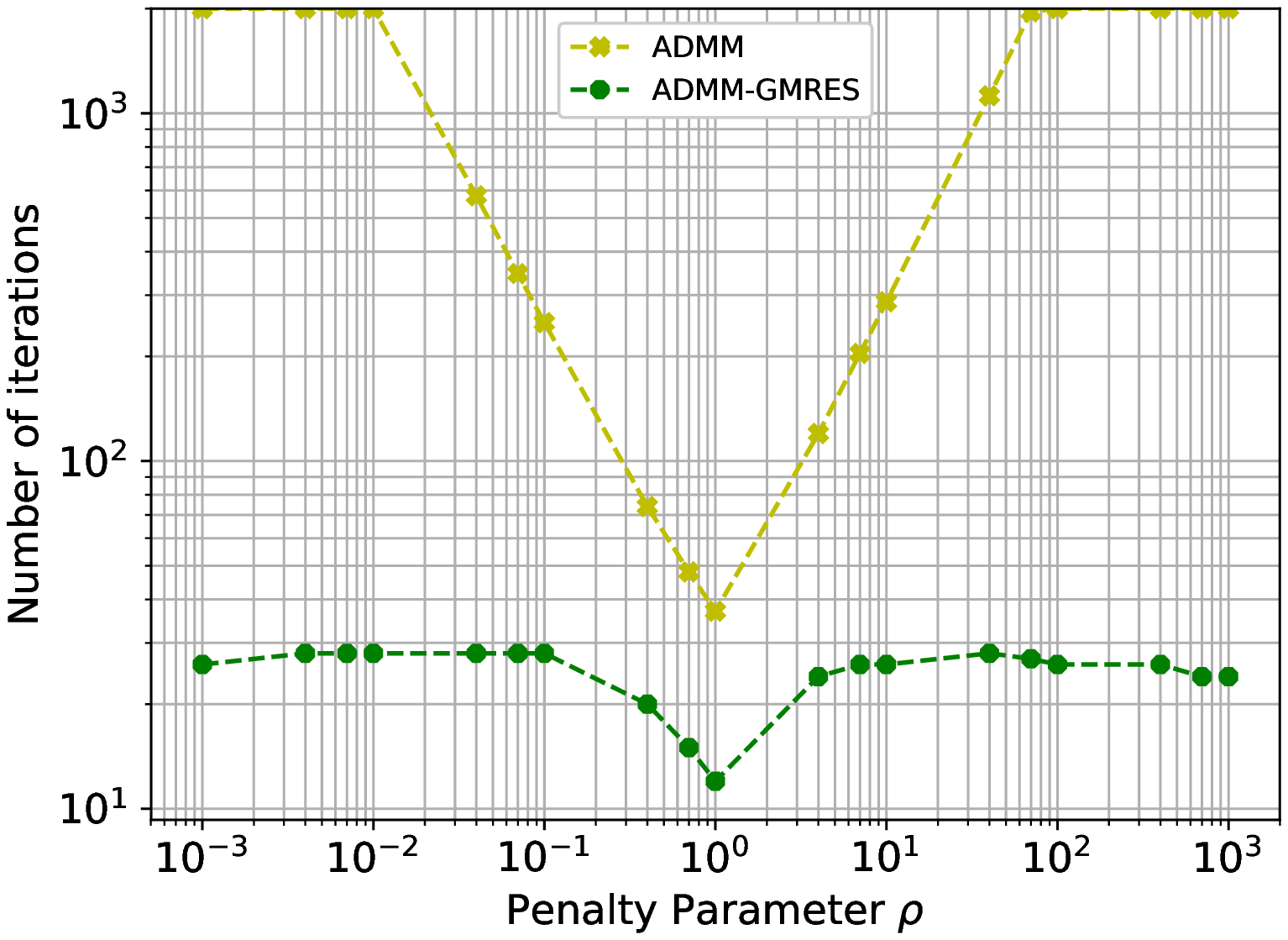}}\label{fig:rho-sens-iterations}}
  \subfigure[KKT system residual]{{\includegraphics[width=3.25in]{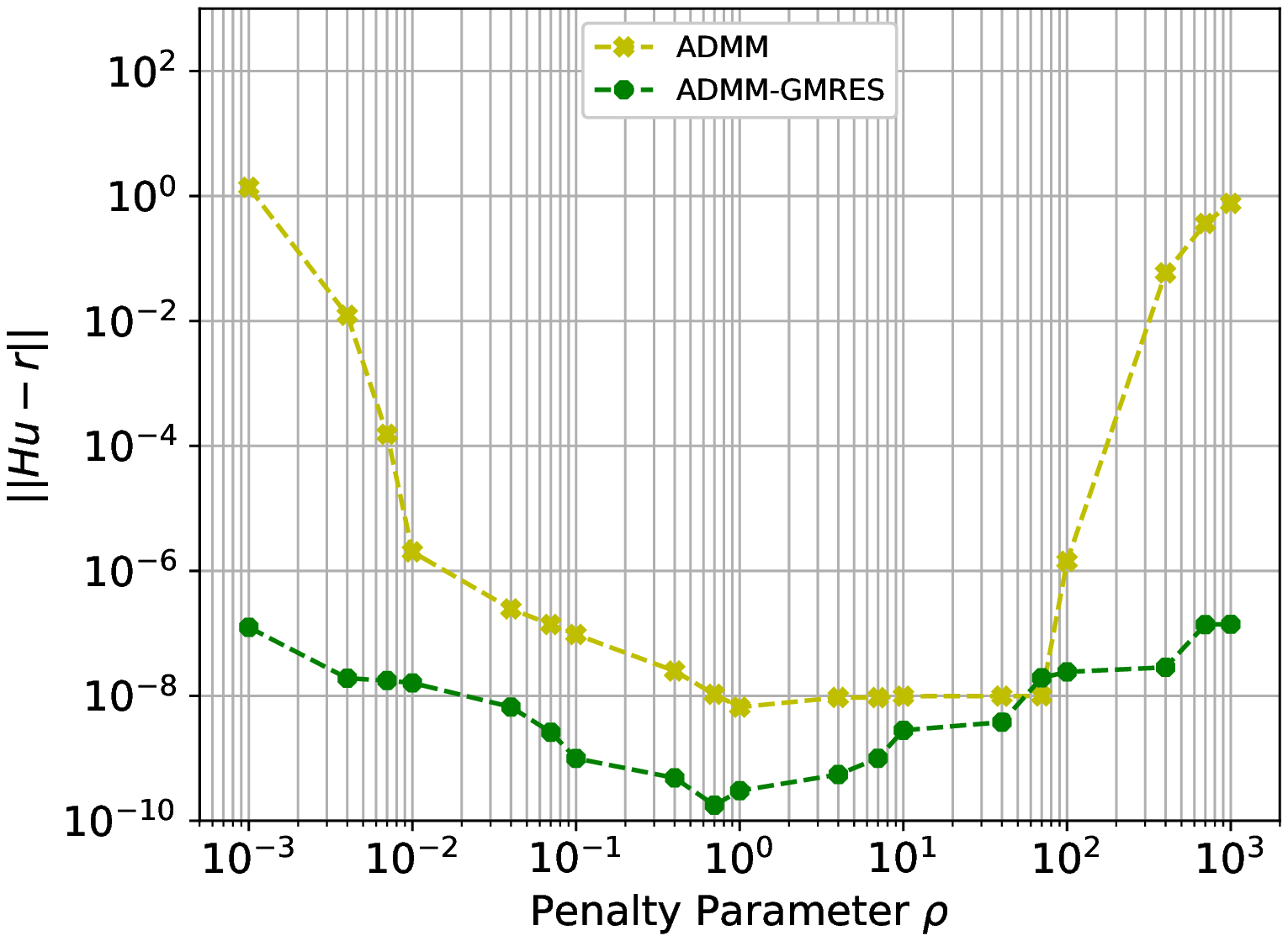}}\label{fig:rho-sense-accuracy}}\\
  \caption{Sensitivity of ADMM and ADMM-GMRES to penalty parameter $\rho$. }
  \label{fig:rho-sens-results}
\end{figure}   

Figure \ref{fig:rho-sens-results} compares the performance of ADMM against that of ADMM-GMRES for a stochastic program with $n_z{=}1{,}000$ and different values of $\rho$ ranging from $10^{-3}$ to $10^{3}$. Here we measure performance in terms of the number of iterations. Our results on the block-structured problem are in agreement with those in \citep{ZhangRichard2016PIiA}, where a single block problem is solved. We see that ADMM-GMRES is remarkably robust to the choice of $\rho$ (the number of iterations remain below 30). ADMM, on the other hand, fails to converge within 2,000 iterations for small and large values of $\rho$. The superior performance of ADMM-GMRES is attributed to the fact that the selection of $\rho$ only has an effect on the preconditioner and not on the convergence properties of GMRES. Nevertheless, we do observe that an optimal selection of $\rho$ improves performance of both ADMM and ADMM-GMRES. The robustness of ADMM-GMRES is a desirable feature when using the solver within more advanced SQP-based solvers. In particular, recent developments of augmented lagrangian interior-point approaches \citep{ChiangNai-Yuan2017AALF, KuhlmannRenke2018ApaL} provide promising frameworks for ADMM-GMRES because the selection of $\rho$ is made based on information from the outer-iteration of the interior-point \citep{ArmandPaul2014Soap, ArmandPaul2017Agaq}. Finally, recent developments of inertia-free methods for nonconvex nonlinear optimization enable the use of iterative linear solvers for the Newton subproblem.

\subsection{Optimal Power Flow Problems}\label{sec:powerflow}

We now demonstrate the computational benefits of using ADMM-GMRES by solving stochastic optimal power flow problems. The optimal power flow problem is frequently used in power networks to determine an efficient dispatch of power generators that satisfy demand and maintains feasible operation conditions. This approach assumes that demand forecasts are accurate and determines a nominal operation point for power generation, power flow in transmission lines and voltage angle at each bus in the power grid. We solve the DC-power flow problem for 35 benchmark problems available in the {\tt PGLIB\_OPF} library \citep{Zimmerman11matpower} and determine nominal operation points for each of them. The solution of each benchmark problem was also obtained using {\tt Ipopt}.

\noindent To asses the computational performance of ADMM-GMRES, we formulated a set-point problem that uses the nominal solution of the DC-power flow problem but seeks to minimize the effect of potential uncertainty in electricity demand values. The optimization solved is the quadratic problem: 
\begin{subequations}
\begin{align}
\underset{x_i, z}{\text{min}}
 \;\; & \;\;   \sum_{s \in \mathcal{P}}\sum_{j \in \Omega_G} w_j(P_{G_{j,s}}-P_{G_j}^{\dag})^2 + \sum_{i,j \in \Omega_L} w_{i,j}(P_{F_{i,j,s}}-P_{F_{i,j}}^{\dag})^2 +  \sum_{j \in \Omega_{\mathcal{B}}} w_j(\theta_{j,s}-\theta_j^{\dag})^2&  \\
\text{s.t.} \;\; &\;\; \sum_{j \in \Omega_{G_i}} P_{G_{j,s}} - P_{L_{i,s}}=\sum_{j \in \Omega_i} P_{F_{i,j,s}}, \;\; i \in \Omega_{\mathcal{B}},  s\in \mathcal{P} \\
&\;\; P_{F_{i,j,s}} =  \frac{\theta_{i,s}-\theta_{j,s}}{X_{i,j}}, \quad \; i,j \in \Omega_L,  s\in \mathcal{P}  \\
&\;\; P_{G_{j,s}}- z_j = 0, \;\quad  j \in \Omega_G ,  s\in \mathcal{P}\label{eq:generator-first-stage}\\
&\;\; \theta_{0,s} = 0, \; \quad s\in \mathcal{P}. 
\end{align}
\label{eq:setpoint}
\end{subequations}

\noindent Here, $\Omega_{\mathcal{B}}$ and $\Omega_L$ denote the set of buses and transmission lines in the network, $\Omega_{G}$ the set of generators,  $\Omega_{G_i}$ the set of generators at bus $i$, and $\Omega_i$ the set of buses connected to bus $i$. The variables in the model are the generator outputs  $P_{G}$ , the flow in the transmission lines $P_F$, and the voltage angles $\theta_j$. As parameters we have the reactance of the lines $X_{i, j}$, the loads $P_{L}$, and the set-point values $P_{G}^{\dag}$,  $P_{F}^{\dag}$, and $\theta^{\dag}$ obtained from the DC-power flow solution. We denote the reference bus as $\theta_0$ and define objective weight values $w$. The goal of formulation \eqref{eq:setpoint} is to find the closest feasible operation to the optimal DC-power flow solution while accounting for potential uncertainty in the demands. In this problem the first-stage variables are the output of the generators for which we use Equation \eqref{eq:generator-first-stage} to enforce the same power generation across the set of scenarios. The dimensionality of the first-stage of this problem is given by the number of generators in the power network. Hence, this number varies from three to $4{,}092$ in the $35$ different benchmark problems considered in our study. For each benchmark we generated $50$ random scenarios with normal random distributed noise on the load $P_{L}$

Tables  \ref{tab:larger-instances} and \ref{tab:smaller-instances} summarize the results for the $35$ benchmarks. The problems were sorted according to their number of coupling variables. Table \ref{tab:larger-instances} presents the results for benchmarks with a first-stage dimension greater than $100$. Results for the smaller benchmarks are shown in Table \ref{tab:smaller-instances}. We see that, for the smaller problems, Schur decomposition is the best alternative as it can give exact solutions in about the same time as ADMM. However, for all benchmarks shown in Table \ref{tab:larger-instances}, ADMM-GMRES finds an $\epsilon$-accurate solution in less time than ADMM and Schur decomposition. In particular, for case13659\_pegase (with $4{,}091$ first-stage variables), ADMM-GMRES solved the problem almost an order of magnitude faster than Schur decomposition. By observing the trend for the rest of the problems, we conclude that these favorable scalability results can be expected to hold as the number of coupling variables increases. We highlight that, for many of the problems shown in Table \ref{tab:larger-instances}, ADMM does not converge after $2{,}000$ iterations; ADMM-GMRES, on the other hand, consistently achieves $\epsilon$-accurate solutions in few iterations and regardless of the choice of $\rho$. In summary, our results demonstrate that ADMM-GMRES provides a plausible approach to overcome the limitations of Schur complement decomposition. 

\begin{sidewaystable}[!htp]
\begin{center}
\setlength{\arrayrulewidth}{1.5pt}
\begin{tabular}{lccccccccccc}
    \hline
    Pglib-matpower Case & $n_x$  & $n_z$ & \multicolumn{3}{c}{Solution Time (s)} & \multicolumn{2}{c}{Iteration Count} & \multicolumn{3}{c}{$\|Hu-r\|$}\\
    \cline{4-6}\cline{6-8}\cline{8-11}
    &&&Schur&ADMM&ADGM&ADMM&ADGM&Schur&ADMM&ADGM\\
    \hline
P1. case13659\_pegase & 2115450 & 4091 & 1670.692 & 723.158 & 183.298 & 2000 & 325 & 2.798E-09 & 3.466E-03 & 1.602E-08 \\ 
P2. case9241\_pegase & 1408950 & 1444 & 369.058 & 507.632 & 81.278 & 2000 & 165 & 3.841E-09 & 3.869E-05 & 2.981E-09 \\ 
P3. case6470\_rte & 849800 & 760 & 124.416 & 336.881 & 51.106 & 2000 & 139 & 2.981E-10 & 3.932E-03 & 1.984E-10 \\ 
P4. case6515\_rte & 845950 & 683 & 112.975 & 333.421 & 40.883 & 2000 & 109 & 3.206E-10 & 7.631E-04 & 2.627E-10 \\ 
P5. case6495\_rte & 843650 & 679 & 111.035 & 320.943 & 37.302 & 2000 & 108 & 2.427E-10 & 7.128E-04 & 4.291E-10 \\ 
P6. case2868\_rte & 389850 & 560 & 49.152 & 174.801 & 23.113 & 2000 & 105 & 3.228E-10 & 4.195E-04 & 1.404E-10 \\ 
P7. case2848\_rte & 382250 & 510 & 43.584 & 168.000 & 22.508 & 2000 & 100 & 1.631E-10 & 2.329E-04 & 4.032E-11 \\ 
P8. case2869\_pegase & 423500 & 509 & 48.424 & 189.774 & 38.293 & 2000 & 154 & 4.435E-10 & 1.135E-03 & 6.830E-10 \\ 
P9. case6468\_rte & 813250 & 398 & 62.343 & 326.375 & 24.449 & 2000 & 65 & 2.774E-10 & 5.404E-04 & 4.416E-10 \\ 
P10. case3012wp\_k & 367050 & 378 & 32.732 & 174.999 & 18.422 & 2000 & 74 & 7.032E-11 & 1.356E-05 & 1.338E-11 \\ 
P11. case1951\_rte & 263900 & 365 & 24.864 & 132.780 & 19.566 & 2000 & 80 & 2.764E-10 & 2.836E-05 & 3.227E-10 \\ 
P12. case2383wp\_k & 295900 & 319 & 24.458 & 152.102 & 17.499 & 2000 & 58 & 6.209E-11 & 7.444E-07 & 1.632E-10 \\ 
P13. case1888\_rte & 249900 & 289 & 19.114 & 129.940 & 16.811 & 2000 & 66 & 1.534E-10 & 2.708E-06 & 2.935E-10 \\ 
P14. case3120sp\_k & 367900 & 272 & 23.729 & 175.397 & 12.998 & 2000 & 56 & 3.402E-11 & 1.214E-08 & 2.220E-10 \\ 
P15. case1354\_pegase & 193200 & 259 & 15.185 & 112.178 & 6.520 & 2000 & 47 & 1.453E-10 & 3.649E-07 & 2.728E-10 \\ 
P16. case240\_pserc & 48650 & 142 & 5.262 & 64.665 & 4.737 & 2000 & 59 & 3.942E-10 & 1.902E-08 & 7.815E-10 \\ 
P17. case2746wp\_k & 311600 & 103 & 8.455 & 100.004 & 11.080 & 1287 & 35 & 4.171E-11 & 9.925E-10 & 4.428E-10 \\
    \hline
\end{tabular}
\end{center}
\caption{Performance for ADMM-GMRES (ADGM), ADMM, and Schur decomposition (Schur) on problem \eqref{eq:setpoint} with $n_z \geq 100$.  \label{tab:larger-instances}}
\end{sidewaystable}

\begin{sidewaystable}[!htp]
\begin{center}
\setlength{\arrayrulewidth}{1.5pt}
\begin{tabular}{lccccccccccc}
    \hline
    Pglib-matpower Case & $n_x$  & $n_z$ & \multicolumn{3}{c}{Solution Time (s)} & \multicolumn{2}{c}{Iteration Count} & \multicolumn{3}{c}{$\|Hu-r\|$}\\
    \cline{4-6}\cline{6-8}\cline{8-11}
    &&&Schur&ADMM&ADGM&ADMM&ADGM&Schur&ADMM&ADGM\\
    \hline
P18. case73\_ieee\_rts & 19200 & 95 & 3.161 & 21.751 & 2.914 & 772 & 38 & 1.001E-11 & 9.795E-10 & 6.237E-11 \\ 
P19. case2746wop\_k & 311100 & 84 & 7.033 & 78.842 & 10.125 & 1022 & 34 & 2.763E-11 & 9.911E-10 & 1.306E-09 \\ 
P20. case2736sp\_k & 308400 & 81 & 6.758 & 81.199 & 9.796 & 1054 & 35 & 2.538E-11 & 9.978E-10 & 4.499E-10 \\ 
P21. case300\_ieee & 41200 & 56 & 2.212 & 41.976 & 6.157 & 1305 & 38 & 6.341E-11 & 9.891E-10 & 6.844E-10 \\ 
P22. case2737sop\_k & 305650 & 53 & 4.668 & 54.118 & 9.791 & 696 & 34 & 1.808E-11 & 9.905E-10 & 3.925E-10 \\ 
P23. case200\_pserc & 26000 & 37 & 1.440 & 13.064 & 3.150 & 428 & 28 & 3.524E-12 & 9.645E-10 & 2.619E-11 \\ 
P24. case24\_ieee\_rts & 6250 & 31 & 1.096 & 6.219 & 2.059 & 231 & 30 & 3.218E-12 & 9.918E-10 & 5.142E-11 \\ 
P25. case118\_ieee & 17050 & 18 & 0.799 & 10.767 & 2.212 & 375 & 32 & 1.329E-11 & 9.662E-10 & 2.874E-10 \\ 
P26. case162\_ieee\_dtc & 23450 & 11 & 0.618 & 4.629 & 1.829 & 148 & 20 & 7.489E-12 & 8.752E-10 & 1.347E-10 \\ 
P27. case89\_pegase & 16100 & 11 & 0.588 & 3.927 & 1.936 & 132 & 23 & 4.015E-12 & 9.344E-10 & 1.228E-10 \\ 
P28. case39\_epri & 5200 & 9 & 0.476 & 3.471 & 1.516 & 124 & 18 & 8.340E-12 & 8.749E-10 & 1.543E-10 \\ 
P29. case30\_as & 4100 & 5 & 0.360 & 2.386 & 0.573 & 81 & 11 & 1.632E-13 & 9.865E-09 & 5.271E-11 \\ 
P30. case30\_fsr & 4100 & 5 & 0.362 & 2.047 & 0.577 & 73 & 11 & 1.997E-13 & 8.138E-09 & 3.175E-11 \\ 
P31. case5\_pjm & 1000 & 4 & 0.319 & 2.707 & 0.474 & 103 & 8 & 2.434E-13 & 9.123E-09 & 3.137E-10 \\ 
P32. case57\_ieee & 7200 & 3 & 0.314 & 2.382 & 0.747 & 82 & 7 & 7.034E-13 & 9.208E-09 & 2.171E-09 \\ 
P33. case14\_ieee & 1850 & 1 & 0.237 & 1.343 & 0.346 & 48 & 3 & 1.355E-13 & 9.339E-09 & 3.045E-11 \\ 
P34. case30\_ieee & 3700 & 1 & 0.243 & 1.357 & 0.354 & 47 & 3 & 5.790E-14 & 9.653E-09 & 6.383E-12 \\ 
P35. case3\_lmbd & 450 & 1 & 0.235 & 1.233 & 0.342 & 44 & 3 & 1.515E-13 & 7.080E-09 & 1.520E-12 \\ 
    \hline
\end{tabular}
\end{center}
\caption{Performance for ADMM-GMRES (ADGM), ADMM and Schur decomposition (Schur) on problem \eqref{eq:setpoint} with $n_z \leq 100$. \label{tab:smaller-instances}}
\end{sidewaystable}

\newpage
%\noindent {\color{red} VICTOR:} The following figures are not referenced in the text (yet). I put them here in case you want to include them in the paper. There is a copy with the y-axis in log scale. I think is better to show this in regular scale since it makes more clear the time differences. However, if you rather present the log scale figures just append the suffix log\_ to the name of the figures. 

\begin{figure}[!htb]
  \begin{center}
   \includegraphics[width=0.7\textwidth]{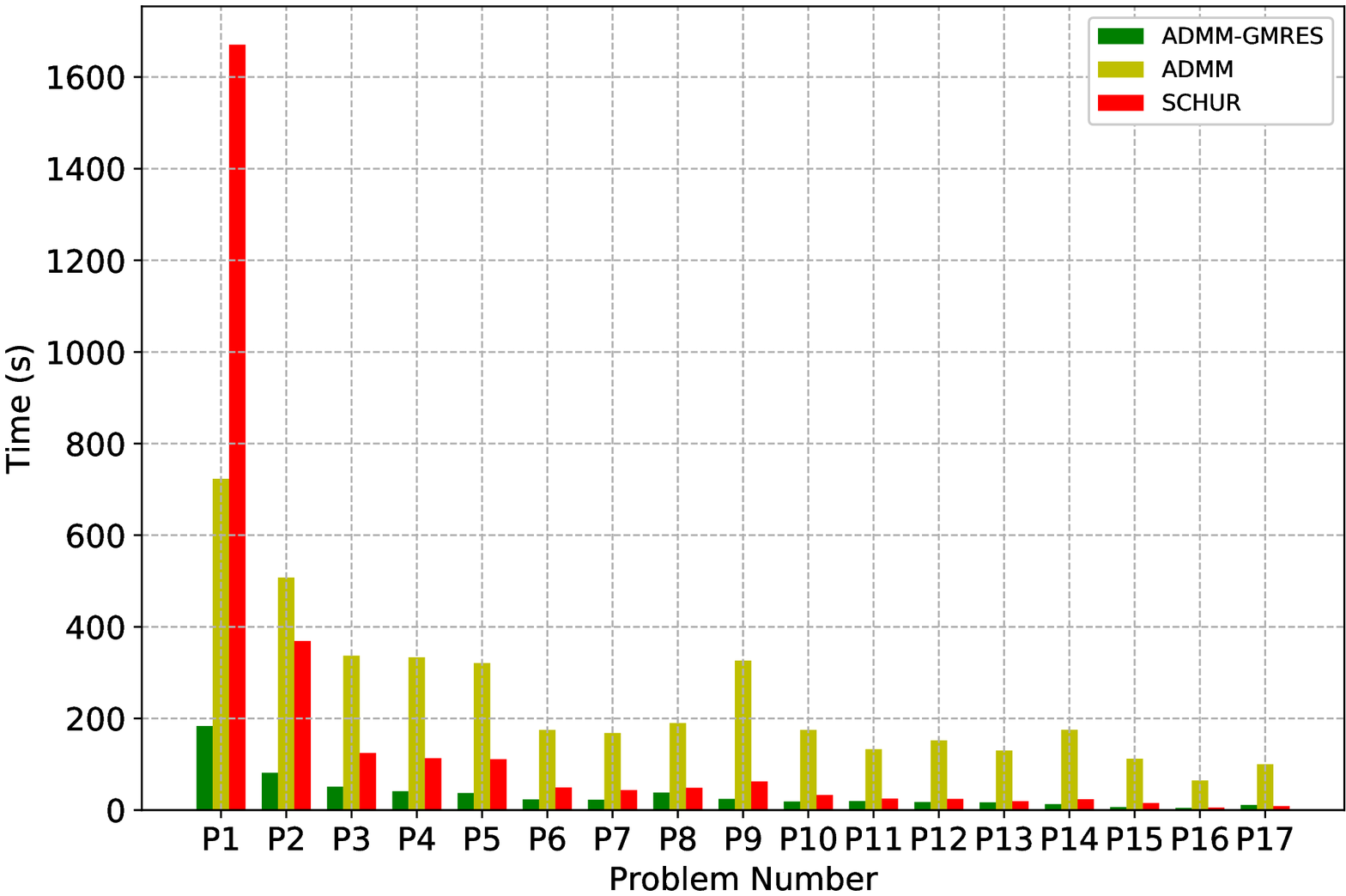}
 \caption{Computational times for Schur decomposition, ADMM, and ADMM-GMRES for problems with $n_z\geq 100$.}\label{fig:times-large-results} 
 \end{center}
\end{figure} 

\begin{figure}[!htb]
  \begin{center}
   \includegraphics[width=0.7\textwidth]{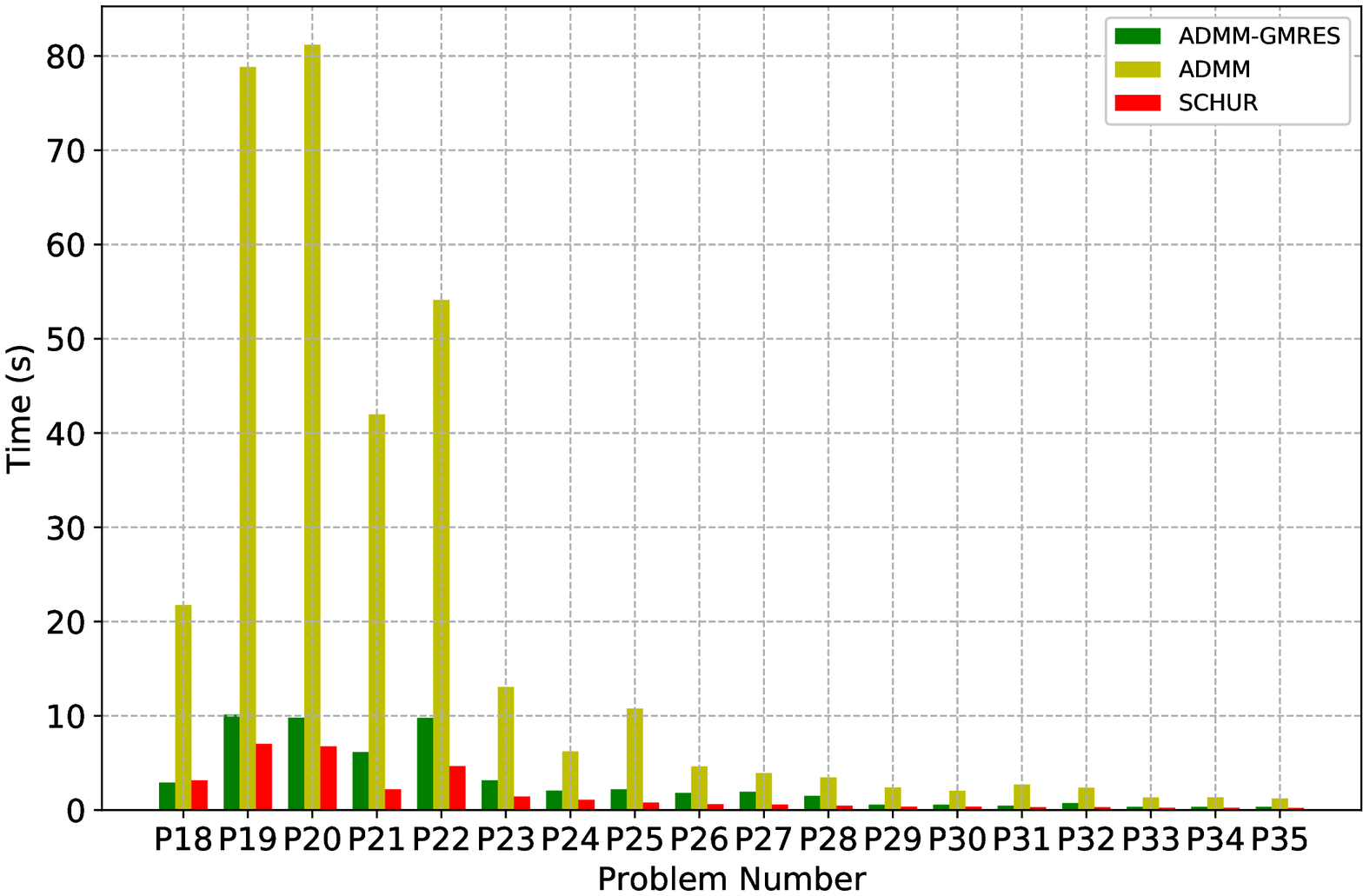}
 \caption{Computational times for Schur decomposition, ADMM, and ADMM-GMRES for problems with $n_z < 100$.}\label{fig:times-small-results} 
 \end{center}
\end{figure} 

\begin{figure}[!htb]
  \begin{center}
   \includegraphics[width=0.7\textwidth]{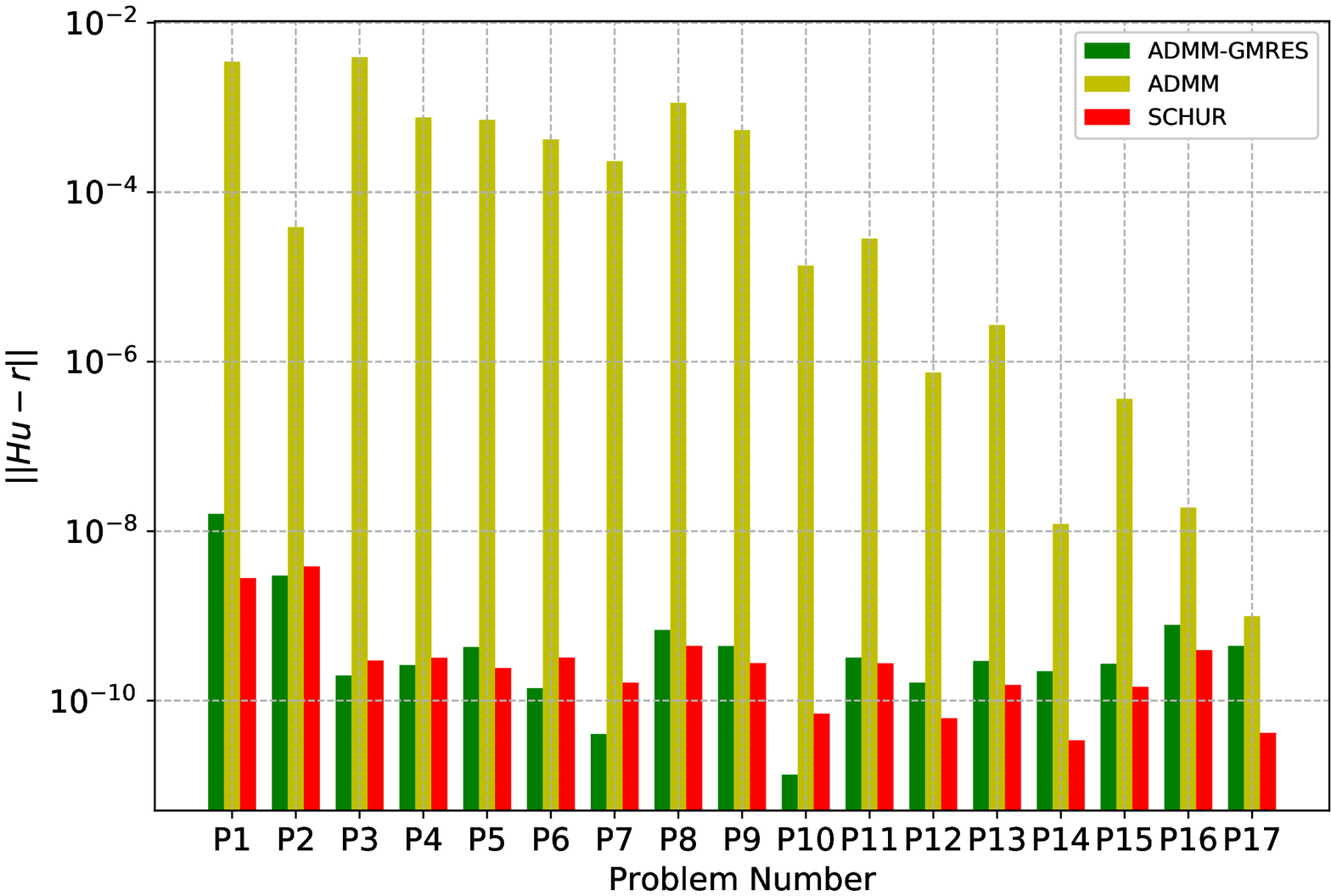}
 \caption{Residuals for Schur decomposition, ADMM, and ADMM-GMRES for problems with $n_z\geq 100$.}\label{fig:residuals-large-results} 
 \end{center}
\end{figure} 

\begin{figure}[!htb]
  \begin{center}
   \includegraphics[width=0.7\textwidth]{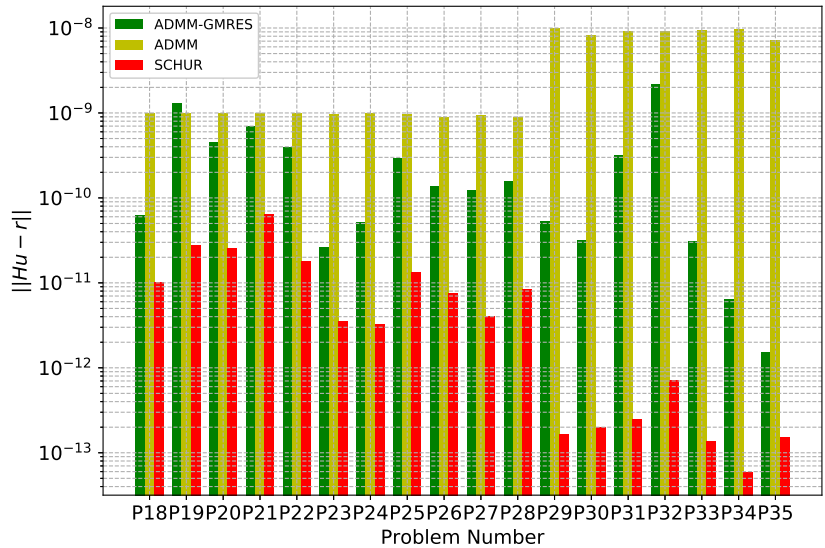}
 \caption{Residuals for Schur decomposition, ADMM, and ADMM-GMRES for problems with $n_z < 100$.}\label{fig:residuals-small-results} 
 \end{center}
\end{figure} 

\section{Conclusions and Future Work}\label{sec:conclusions}

We have demonstrated that ADMM provides an effective mechanism to precondition iterative linear solvers and with this overcome scalability limitations of Schur complement decomposition.  Our results also demonstrate that the approach is robust to the choice of the penalty parameter. As part of future work, we will investigate the performance of ADMM-GMRES within a nonlinear interior-point framework. Here, it will be necessary to relax our assumptions on strong convexity and on the full rank of the Jacobian. Preliminary results reported  in the literature indicate that different types of primal-dual regularized KKT systems can be used to compute search steps within interior-point methods under such relaxed conditions \citep{ChiangNai-Yuan2017AALF}. For instance, the primal-dual regularized system correspond to the optimality conditions of the QP problem:

\begin{equation} 
\label{eq:ecqp-regularized}
\begin{array}{rrclcl} 
\displaystyle \min_{x, z, r} & \multicolumn{3}{l}{\frac{1}{2}x^T (D +\delta I) x + c^T x + \frac{1}{2\rho}\|r\|^2 + \frac{\rho}{2}\|A x  + B z - \frac{1}{\rho}r\|^2 }\\
\textrm{s.t.} & A x  + B z - \frac{1}{\rho}r& = & 0, \;\;\; (y)
\end{array} 
\end{equation} 

We will investigate ADMM variants to precondition such systems. The effectiveness of using ADMM as a preconditioner makes us wonder whether other approaches can be used for preconditioning as well. For instance, inexact dual Newton strategies can potentially be used to precondition structured KKT systems. This is an interesting direction of future work. We will also investigate advanced ADMM strategies that use second-order multiplier updates to accelerate the preconditioner. 

\section*{Disclaimer}\label{sec:disclaimer}
Sandia National Laboratories is a multimission laboratory managed and operated by National Technology and Engineering Solutions of Sandia, LLC., a wholly owned subsidiary of Honeywell International, Inc., for the U.S. Department of Energy's National Nuclear Security Administration under contract  {DE-NA-0003525}. This paper describes objective technical results and analysis. Any subjective views or opinions that might be expressed in the paper do not necessarily represent the views of the U.S. Department of Energy or the United States Government.This work was conducted as part of the Institute for the Design of Advanced Energy Systems (IDAES) with funding from the Office of Fossil Energy, Cross-Cutting Research, U.S. Department of Energy. V. M. Zavala acknowledges funding from the National Science Foundation under award NSF- EECS-1609183.

\bibliography{references}

\begin{thebibliography}{10}

\bibitem{MUMPS1}
P.~R. Amestoy, I.~S. Duff, J.~Koster, and J.-Y. L'Excellent.
\newblock A fully asynchronous multifrontal solver using distributed dynamic
  scheduling.
\newblock {\em SIAM Journal on Matrix Analysis and Applications}, 23(1):15--41,
  2001.

\bibitem{MUMPS2}
P.~R. Amestoy, A.~Guermouche, J.-Y. L'Excellent, and S.~Pralet.
\newblock Hybrid scheduling for the parallel solution of linear systems.
\newblock {\em Parallel Computing}, 32(2):136--156, 2006.

\bibitem{ArmandPaul2014Soap}
Paul Armand, Joël Benoist, Riadh Omheni, and Vincent Pateloup.
\newblock Study of a primal-dual algorithm for equality constrained
  minimization.
\newblock {\em Computational Optimization and Applications}, 59(3):405--433,
  2014.

\bibitem{ArmandPaul2017Agaq}
Paul Armand and Riadh Omheni.
\newblock A globally and quadratically convergent primal–dual augmented
  lagrangian algorithm for equality constrained optimization.
\newblock {\em Optimization Methods and Software}, 32(1):1--21, 2017.

\bibitem{BenziMichele2005Nsos}
Michele Benzi, Gene~H. Golub, and Jrg Liesen.
\newblock Numerical solution of saddle point problems.
\newblock {\em Acta Numerica}, 14:1--137, 2005.

\bibitem{BenziMichele2006Oteo}
Michele Benzi and Valeria Simoncini.
\newblock On the eigenvalues of a class of saddle point matrices.
\newblock {\em Numerische Mathematik}, 103(2):173--196, 2006.

\bibitem{boyd2011distributed}
Stephen Boyd, Neal Parikh, Eric Chu, Borja Peleato, Jonathan Eckstein, et~al.
\newblock Distributed optimization and statistical learning via the alternating
  direction method of multipliers.
\newblock {\em Foundations and Trends{\textregistered} in Machine learning},
  3(1):1--122, 2011.

\bibitem{CaoYankai2016Cpfs}
Yankai Cao, Carl Laird, and Victor Zavala.
\newblock Clustering-based preconditioning for stochastic programs.
\newblock {\em Computational Optimization and Applications}, 64(2):379--406,
  2016.

\bibitem{ChiangNai-Yuan2017AALF}
Nai-Yuan Chiang, Rui Huang, and Victor~M. Zavala.
\newblock An augmented lagrangian filter method for real-time embedded
  optimization.
\newblock {\em Automatic Control, IEEE Transactions on}, 62(12):6110--6121,
  2017.

\bibitem{pipsnlp}
Naiyuan Chiang, Cosmin~G Petra, and Victor~M Zavala.
\newblock Structured nonconvex optimization of large-scale energy systems using
  {PIPS-NLP}.
\newblock In {\em Proc. of the 18th Power Systems Computation Conference
  (PSCC), Wroclaw, Poland}, 2014.

\bibitem{eckstein1992douglas}
Jonathan Eckstein and Dimitri~P Bertsekas.
\newblock On the {Douglas-Rachford} splitting method and the proximal point
  algorithm for maximal monotone operators.
\newblock {\em Mathematical Programming}, 55(1-3):293--318, 1992.

\bibitem{ElmanHowardC.1994IaPU}
Howard~C. Elman and Gene~H. Golub.
\newblock Inexact and preconditioned uzawa algorithms for saddle point
  problems.
\newblock {\em SIAM Journal on Numerical Analysis}, 31(6):1645--1661, 1994.

\bibitem{ForsgrenAnders2007ISoA}
Anders Forsgren, Philip~E. Gill, and Joshua~D. Griffin.
\newblock Iterative solution of augmented systems arising in interior methods.
\newblock {\em SIAM Journal on Optimization}, 18(2):666--690, 2007.

\bibitem{ampl}
R.~Fourer, D.M. Gay, and B.W. Kernighan.
\newblock {\em AMPL: A Modeling Language for Mathematical Programming}.
\newblock Scientific Press, 1993.

\bibitem{GaulAndre2015PRKs}
André Gaul and Nico Schlömer.
\newblock Preconditioned recycling krylov subspace methods for self-adjoint
  problems.
\newblock {\em arXiv.org}, 2015.

\bibitem{GhadimiEuhanna2015OPSf}
Euhanna Ghadimi, Andre Teixeira, Iman Shames, and Mikael Johansson.
\newblock Optimal parameter selection for the alternating direction method of
  multipliers (admm): Quadratic problems.
\newblock {\em IEEE Transactions on Automatic Control}, 60(3):644--658, 2015.

\bibitem{GillPhilipE.1992PfIS}
Philip~E. Gill, Walter Murray, Dulce~B. Poncelen, and Michael~A. Saunders.
\newblock Preconditioners for indefinite systems arising in optimization.
\newblock {\em SIAM Journal on Matrix Analysis and Applications},
  13(1):292--311, 1992.

\bibitem{GoldfarbDonald2012FMAf}
Donald Goldfarb and Shiqian Ma.
\newblock Fast multiple-splitting algorithms for convex optimization.
\newblock {\em SIAM Journal on Optimization}, 22(2):533--556, 2012.

\bibitem{GolubGeneH.2003OSBI}
Gene~H. Golub and Chen Greif.
\newblock On solving block-structured indefinite linear systems.
\newblock {\em SIAM Journal on Scientific Computing}, 24(6):2076--2092, 2003.

\bibitem{GondzioJacek2009Esip}
J.~Gondzio and A.~Grothey.
\newblock Exploiting structure in parallel implementation of interior point
  methods for optimization.
\newblock {\em Computational Management Science}, 6(2):135--160, May 2009.

\bibitem{GondzioJacek2003Pisf}
J.~Gondzio and R.~Sarkissian.
\newblock Parallel interior-point solver for structured linear programs.
\newblock {\em Mathematical Programming}, 96(3):561--584, June 2003.

\bibitem{GuoKe2017CoAf}
Ke~Guo, Deren Han, David Wang, and Tingting Wu.
\newblock Convergence of admm for multi-block nonconvex separable optimization
  models.
\newblock {\em Frontiers of Mathematics in China}, 12(5):1139--1162, 2017.

\bibitem{HanDeren2013LLCo}
Deren Han and Xiaoming Yuan.
\newblock Local linear convergence of the alternating direction method of
  multipliers for quadratic programs.
\newblock {\em SIAM Journal on Numerical Analysis}, 51(6):3446--3457, 2013.

\bibitem{hart2017pyomo}
W.~E. Hart, C.~D. Laird, J.~Watson, D.~L. Woodruff, G.~A. Hackebeil, B.~L.
  Nicholson, and J.~D. Siirola.
\newblock {\em Pyomo--optimization modeling in python}, volume~67.
\newblock Springer Science and Business Media, second edition, 2017.

\bibitem{HeBingsheng2012OtOC}
Bingsheng He and Xiaoming Yuan.
\newblock On the $o(1/n)$ convergence rate of the douglasrachford alternating
  direction method.
\newblock {\em SIAM Journal on Numerical Analysis}, 50(2):700--709, 2012.

\bibitem{HongMingyi2017Otlc}
Mingyi Hong and Zhi-Quan Luo.
\newblock On the linear convergence of the alternating direction method of
  multipliers.
\newblock {\em Mathematical Programming}, 162(1):165--199, March 2017.

\bibitem{HSL}
Eric Jones, Travis Oliphant, Pearu Peterson, et~al.
\newblock A collection of fortran codes for large scale scientific computation.
\newblock [Online; accessed 11/03/2019].

\bibitem{scipyref}
Eric Jones, Travis Oliphant, Pearu Peterson, et~al.
\newblock {SciPy}: Open source scientific tools for {Python}, 2001--.
\newblock [Online; accessed 11/03/2019].

\bibitem{Kang2014}
Jia Kang, Yankai Cao, Daniel~P. Word, and C.~D. Laird.
\newblock {An interior-point method for efficient solution of block-structured
  {NLP} problems using an implicit Schur-complement decomposition}.
\newblock {\em Computers and Chemical Engineering}, 71:563--573, 2014.

\bibitem{KuhlmannRenke2018ApaL}
Renke Kuhlmann and Christof Büskens.
\newblock A primal–dual augmented lagrangian penalty-interior-point filter
  line search algorithm.
\newblock {\em Mathematical Methods of Operations Research}, 87(3):451--483,
  2018.

\bibitem{MaC.F.2015TcUm}
C.~F. Ma and Q.~Q. Zheng.
\newblock The corrected uzawa method for solving saddle point problems.
\newblock {\em Numerical Linear Algebra with Applications}, 22(4):717--730,
  2015.

\bibitem{MoralesJoseLuis2000APbL}
Jose~Luis Morales and Jorge Nocedal.
\newblock Automatic preconditioning by limited memory quasi-newton updating.
\newblock {\em SIAM Journal on Optimization}, 10(4), 2000.

\bibitem{QuarteroniAlfio2007Nm}
Alfio Quarteroni.
\newblock {\em Numerical mathematics}.
\newblock Springer, Berlin ; New York, 2nd ed.. edition, 2007.

\bibitem{RodriguezJoseS.2018BAin}
Jose~S. Rodriguez, Bethany Nicholson, Carl Laird, and Victor~M. Zavala.
\newblock Benchmarking admm in nonconvex nlps.
\newblock {\em Computers and Chemical Engineering}, 119:315--325, 2018.

\bibitem{RustenTorgeir1992APIM}
Torgeir Rusten and Ragnar Winther.
\newblock A preconditioned iterative method for saddlepoint problems.
\newblock {\em SIAM Journal on Matrix Analysis and Applications},
  13(3):887--904, 1992.

\bibitem{BrendtWohlberg2017APPS}
Brendt Wohlberg.
\newblock Admm penalty parameter selection by residual balancing.
\newblock {\em arXiv.org}, 2017.

\bibitem{WordDaniel2014Epso}
D.~P. Word.
\newblock Efficient parallel solution of large-scale nonlinear dynamic
  optimization problems.
\newblock {\em Computational Optimization and Applications}, 59(3):667--689,
  December 2014.

\bibitem{zavalalaird}
Victor~M Zavala, Carl~D Laird, and Lorenz~T Biegler.
\newblock Interior-point decomposition approaches for parallel solution of
  large-scale nonlinear parameter estimation problems.
\newblock {\em Chemical Engineering Science}, 63(19):4834--4845, 2008.

\bibitem{schurbook}
F.~Zhang.
\newblock {\em The Schur Complement and Its Applications}.
\newblock Numerical Methods and Algorithms, 4. Springer US, 2005.

\bibitem{ZhangRichard2016PIiA}
Richard Zhang and Jacob White.
\newblock Parameter insensitivity in admm-preconditioned solution of
  saddle-point problems.
\newblock {\em arXiv.org}, 2016.

\bibitem{ZhangRichardY.2018GAfQ}
Richard~Y. Zhang and Jacob~K. White.
\newblock Gmres-accelerated admm for quadratic objectives.
\newblock {\em SIAM Journal on Optimization}, 28(4):3025--3056, 2018.

\bibitem{Zimmerman11matpower}
Ray~Daniel Zimmerman, Carlos~Edmundo Murillo-sánchez, Robert~John Thomas, and
  Life Fellow.
\newblock Matpower steady-state operations, planning and analysis tools for
  power systems research and education.
\newblock {\em IEEE Transactions on Power Systems}, pages 12--19, 2011.

\bibitem{WalterZulehner2002Aoim}
Walter Zulehner.
\newblock Analysis of iterative methods for saddle point problems: a unified
  approach.
\newblock {\em Mathematics of Computation}, 71(238):479--505, 2002.

\end{thebibliography}

\appendix

\section{Computing Operator $T_\rho(u)$ using ADMM}\label{app:admm-richardson}

Here we prove that the operator $T_\rho(u)$ can be computed by applying one ADMM iteration. Consider, without loss of generality (and in order to simplify the presentation), the case of a single block problem of the form:
\begin{subequations}
\begin{align}
\underset{ x, z}{\text{min}}
 \;\; & \;\;  \frac{1}{2}x^T D x + c^T x + \frac{\rho}{2}\|Ax + Bz\|^{2}& \\
\text{s.t.} \;\; &\;\; Ax + Bz = 0,\; (y).
\end{align}
\label{eq:single-block-qp-augmented}
\end{subequations}
The results that we derive next can be extended to multiple blocks using induction. The KKT system for problem \eqref{eq:single-block-qp-augmented} is:

\begin{equation}
\label{eq:kkt-aug-compact-appendix2}
\begin{bmatrix}
K_{\rho} &   \rho A^{T}B &  A^{T}   \\
 \rho B^TA &       \rho B^{T}B         & B^{T}                \\
A  &   B      &                 \\   
\end{bmatrix}
\begin{bmatrix}
x\\
z\\
y\\
\end{bmatrix}
=
\begin{bmatrix}
-c \\
0 \\
0 \\
\end{bmatrix}
\end{equation}

\noindent where $K_{\rho}=D+\rho A^{T}A$. Applying a Gauss-Seidel splitting to this system at a point $u=(x,z,y)$ leads to the update $u^{+}=T_{\rho}(u)=G_{\rho}u + f_{\rho}$, where $G_{\rho} = M_{\rho}^{-1}N_\rho$ and $f_\rho=M_{\rho}^{-1}r$. The explicit form of  $M_{\rho}^{-1}$ is given by:

\begin{equation}
M_{\rho}^{-1} = 
\begin{bmatrix}
K_{\rho}^{-1} & 0& 0\\
-\Sigma^{-1}B^{T}AK_{\rho}^{-1} & \frac{1}{\rho}\Sigma^{-1} & 0\\
\rho (I-B\Sigma^{-1}B^{T})AK_{\rho}^{-1} & B\Sigma^{-1} & -\rho I
\end{bmatrix}
\end{equation} 

\noindent where $\Sigma:=B^{T}B$. Having $M_{\rho}^{-1}$ we construct:

\begin{align}
G_{\rho} = M_{\rho}^{-1}N_\rho &= 
\begin{bmatrix}
K_{\rho}^{-1} & 0& 0\\
-\Sigma^{-1}B^{T}AK_{\rho}^{-1} & \frac{1}{\rho}\Sigma^{-1} & 0\\
\rho (I-B\Sigma^{-1}B^{T})AK_{\rho}^{-1} & B\Sigma^{-1} & -\rho I
\end{bmatrix}
\begin{bmatrix}
0 & -\rho A^{T}B & -A^{T} \\
0 & 0 & -B^{T} \\ 
0 & 0 & -\frac{1}{\rho}I
\end{bmatrix}\nonumber\\
&= \begin{bmatrix}
0 & -\rho K_{\rho}^{-1}A^{T}B& -K_{\rho}^{-1}A^{T}\\
0 & \Sigma^{-1}B^{T}AK_{\rho}^{-1}A^{T}B & \;\;\Sigma^{-1}B^{T}(AK_{\rho}^{-1}A^{T} - \frac{1}{\rho}I)\\
0 & \rho^2(B\Sigma^{-1}B^{T}-I)AK_{\rho}^{-1}A^{T}B & \rho(B\Sigma^{-1}B^{T}-I)AK_{\rho}^{-1}A^{T} - B\Sigma^{-1}B^{T} + I
\end{bmatrix}
\end{align}

\noindent and the right-hand-side-vector

\begin{align}
f_{\rho} = M_{\rho}^{-1}r &= 
\begin{bmatrix}
K_{\rho}^{-1} & 0& 0\\
-\Sigma^{-1}B^{T}AK_{\rho}^{-1} & \frac{1}{\rho} \Sigma^{-1} & 0\\
\rho (I-B\Sigma^{-1}B^{T})AK_{\rho}^{-1} & B\Sigma^{-1} & -\rho I
\end{bmatrix}
\begin{bmatrix}
-c \\
0  \\ 
0 \\
\end{bmatrix}\qquad\qquad\qquad\qquad\qquad\qquad\quad\nonumber\\
&= \begin{bmatrix}
-K_{\rho}^{-1}c \\
\Sigma^{-1}B^{T}AK_{\rho}^{-1}c \\
\rho (B\Sigma^{-1}B^{T}-I)AK_{\rho}^{-1}c
\end{bmatrix}
\end{align}

\noindent By defining $Q:=AK_{\rho}^{-1}A^{T}$ we can write the explicit form of the update $u^+$ as:

\begin{equation}
\begingroup\makeatletter\def\f@size{10}\check@mathfonts
\begin{bmatrix}
x \\
z  \\ 
y \\
\end{bmatrix}^{+}=
\begin{bmatrix}
0 & -\rho K_{\rho}^{-1}A^{T}B& -K_{\rho}^{-1}A^{T}\\
0 & \Sigma^{-1}B^{T}QB & \;\;\Sigma^{-1}B^{T}(Q - \frac{1}{\rho}I)\\
0 & \rho^2(B\Sigma^{-1}B^{T}-I)QB & \rho (B\Sigma^{-1}B^{T}-I)Q - B\Sigma^{-1}B^{T} + I
\end{bmatrix}
\begin{bmatrix}
x \\
z  \\ 
y \\
\end{bmatrix}
+
\begin{bmatrix}
-K_{\rho}^{-1}c \\
\Sigma^{-1}B^{T}AK_{\rho}^{-1}c \\
\rho (B\Sigma^{-1}B^{T}-I)AK_{\rho}^{-1}c
\end{bmatrix}
\endgroup
\end{equation}

\noindent Upon expansion we obtain the update $u^+=(x^+,z^+,y^+)$:
\begin{subequations}
\begin{align}
x^{+} &= -\rho K_{\rho}^{-1}A^{T}Bz - K_{\rho}^{-1}A^{T}y - K_{\rho}^{-1}c \label{eq:update-x}\\
z^{+} &= \Sigma^{-1}B^{T}QBz + \Sigma^{-1}B^{T}(Q - \frac{1}{\rho}I)y + \Sigma^{-1}B^{T}AK_{\rho}^{-1}c \label{eq:update-z}\\
y^{+} &= \rho^2(B\Sigma^{-1}B^{T}-I)QBz + \left[\rho (B\Sigma^{-1}B^{T}-I)Q - B\Sigma^{-1}B^{T} + I\right]y + \rho (B\Sigma^{-1}B^{T}-I)AK_{\rho}^{-1}c \label{eq:update-y}
\end{align}
\end{subequations}

\noindent We now show that ADMM delivers the same updates after one iteration. We use the augmented Lagrange function:

\begin{equation}
\mathcal{L}(x,z,y) = x^TDx + c^Tx  +\left(Ax  + Bz\right)^{T}y + \frac{\rho}{2}\|Ax  + Bz\|^{2}.
\end{equation}

\noindent Initializing at $u=(x,z,y)$, the update $x^+$ is given by:

\begin{equation}
x^{+} = \argmin\limits_{x}  \mathcal{L}_{\rho}(x,z,y)
\end{equation}
\noindent For which the optimality conditions are
\begin{align}
\nabla_x \mathcal{L}_{\rho}(x,z,y)=(D+\rho A^TA)x+\rho A^{T}Bz+A^{T}y+c=0
\end{align}
\noindent and thus,
\begin{align}
x^{+}&=-(D+ \rho A^TA)^{-1}\left[\rho A^{T}Bz + A^{T}y +c\right]\nonumber\\
&=-K_{\rho}^{-1}\left[\rho A^{T}Bz+A^{T}y+c\right]\nonumber\\
&=-\rho K_{\rho}^{-1}A^{T}Bz-K_{\rho}^{-1}A^{T}y-K_{\rho}^{-1}c\label{eq:admm-updatex}
\end{align}
\noindent We note that \eqref{eq:admm-updatex} and \eqref{eq:update-x} are equivalent. The update for the coupling variables $z^+$ is given by:
\begin{equation}
z^{+} = \argmin\limits_{z}  \mathcal{L}_{\rho}(x^{+},z,y)
\end{equation}
\noindent The optimality conditions are given by:
\begin{align}
\nabla_z \mathcal{L}_{\rho}(x^{+},z,y)=B^{T}y + \rho B^{T}Ax^{+}+\rho B^{T}Bz=0
\end{align}
\noindent and thus,
\begin{align}
z^{+}&=-\left(B^{T}B\right)^{-1}\left[\frac{1}{\rho}B^{T}y+B^{T}Ax^{+}\right]\nonumber\\
&=-\Sigma^{-1}\left[\frac{1}{\rho}B^{T}y+B^{T}Ax^{+}\right]\nonumber\\
&=-\Sigma^{-1}\left[\frac{1}{\rho}B^{T}y-B^{T}A^{-1}A^{T}Bz-B^{T}A^{-1}A^{T}y-B^{T}A^{-1}c\right]\nonumber\\
&=-\Sigma^{-1}\left[\frac{1}{\rho}B^{T}y-B^{T}QBz-B^{T}Qy-B^{T}A^{-1}c)\right]\nonumber\\
&=-\Sigma^{-1}\left[B^{T}(\frac{1}{\rho}I-Q)y-B^{T}QBz-B^{T}A^{-1}c\right]\nonumber\\
&=\Sigma^{-1}B^{T}QBz+\Sigma^{-1}B^{T}(Q-\frac{1}{\rho}I)y+\Sigma^{-1}B^{T}A^{-1}c\label{eq:admm-updatez}
\end{align}
\noindent We thus have that \eqref{eq:admm-updatez} and \eqref{eq:update-z} are equivalent. Finally, the dual variables are updated as $y^{+} = y + \rho (Ax^{+}+Bz^{+})$. Substituting \eqref{eq:admm-updatex} and \eqref{eq:admm-updatez} in this expression leads to \eqref{eq:update-y} .

\end{document}